\documentclass[]{article}
\usepackage{graphicx}
\usepackage{amssymb}
\usepackage{eufrak,mathrsfs,amsmath,longtable}

\newtheorem{thm}{Theorem}

\begin{document}

\title{
	Lorentz Ricci solitons of 4-dimensional non-Abelian nilpotent Lie groups}
\author{Rohollah Bakhshandeh-Chamazkoti,\\[3mm]
{\small Department of Mathematics, Faculty of Basic Sciences, Babol  Noshirvani University of Technology, Babol, Iran.}\\[3mm]
{\small Email: r\_bakhshandeh@nit.ac.ir.}
}
\maketitle
\begin{abstract}
The goal of  this paper is to investigate which one of thenon-isometric  left invariant Lorentz metrics
on 4-dimensional nilpotent Lie groups $H_3 \times {\Bbb R}$ and $G_4$ satisfy in Ricci Soliton equation.  
 Among the left-invariant Lorentzian metrics on $H_3 \times {\Bbb R}$, ~ $g_{\lambda}^{+}$  is a shrinking  
 while  $g_{\lambda}^{+}$  and  $g_{\mu}$ are expanding and also $g_0^1, g_0^2, g_0^3$
 have Ricci solitons. 
 We exhibit among the non-isometric left invariant Lorentz metric on the group $G_4$ only $g_1^\lambda, g_2^\lambda$ have 
 Lorentz Ricci solitons and $g_2^\lambda$ is a shrinking.  
%
\end{abstract}
\vspace{0.7cm}
\section{Introduction}
Rahmani, \cite{Rahmani1992}, showed the Heisenberg group $H_3$ has three non-isometric 
left-invariant Lorentzian metrics $g_1$, $g_2$, and $g_3$ and in, \cite{Rahmani2006}, he and his collaborator presented that the associated Lie algebras
of infinitesimal isometries of $(H_3, g_1)$ and $(H_3, g_2)$ are four-dimensional and solvable but
not nilpotent, and the associated Lie algebra of infinitesimal isometries of $(H_3, g_3)$ is six-dimensional and 
 the left-invariant Lorentzian metric $g_2$ has negative
constant curvature $-\frac{1}{4}$, $g_3$ is flat, and $g_1$ is not Einstein. 
 In, \cite{Onda2010}, Onda characterized the left-invariant Lorentzian metric $g_1$ as a Lorentz Ricci
soliton. Moreover, he proved that the group rigid motions of Euclidean $2$-space, $E(2)$, and the group
 rigid motions of Minkowski  $2$-space, $E(1, 1)$, have Lorentz Ricci solitons.
 
 In, \cite{Wears2016}, Wears carried out a complete investigation of the left invariant Lorentzian
metrics on a five-dimensional, connected, simply-connected, two-step nilpotent Lie group and its left invariant
Ricci soliton and algebraic Ricci soliton metrics. One of his student in her dissertation, \cite{Walker2017}, tried to presented 
classification of Lorentzian scalar products of some four dimensional Lie algebras as  Ricci solitons. 

Recently, Calvaruso, \cite{Calvaruso2021},  studied the semi-direct extensions $G_S =H \rtimes \exp({\Bbb R}S)$ 
of the three-dimensional Heisenberg group $H$, equipped with a one-parameter
family of left-invariant metrics $g_a$, $a^2\neq1$ and he  calculated several curvature properties of $(G_S, g_a)$ and
presented a complete classification of its algebraic Ricci solitons with constructing some new examples of
non-algebraic Lorentzian Ricci solitons.

Bokan and et. al., \cite{Bokan2015}, classified left invariant Lorentz metrics
on 4-dimensional nilpotent Lie groups $H_3 \times {\Bbb R}$ and $G_4$.  
Magnin, \cite{Magnin1986}, proved that, up to isomorphism, $\mathfrak{h}_3 \oplus {\Bbb R}$ and $\mathfrak{g}_4$ 
with corresponding Lie groups $H_3 \times {\Bbb R}$ and $G_4$ are only two non-Abelian nilpotent Lie algebras of dimension 4.

Let $g_0$ be a pseudo-Riemannian metric on manifold $M^n$.
If $g_0$ satisfies
$$2{\rm Ric} [g_0] + {\cal L}_{\bf X} g_0 + \alpha g_0 = 0 \ ,$$
where ${\bf X}$ is some vector field and $\alpha$ is some constant,
then $(M^n, g_0, X, \alpha )$ is called a {\it Ricci soliton structure} and $g_0$ is called {\it the Ricci soliton.}
Moreover we say that the Ricci soliton $g_0$ is a {\it gradient Ricci soliton}
if the vector field ${\bf X}$ satisfies ${\bf X}=\nabla f$, where $f$ is some function,
and the Ricci soliton $g$ is a {\it non-gradient Ricci soliton}
if the vector field ${\bf X}$ satisfies $X\ne \nabla f$ for any function $f$.
If a constant $\alpha$ is negative, zero, or positive,
then $g$ is called a shrinking, steady, or expanding Ricci soliton, respectively.
Recall that the Ricci solitons have a relation with the Ricci flow.
 A pseudo-Riemannian metric $g_0$ is a Ricci soliton if and only if
$g_0$ is a solution of the Ricci flow equation,
\begin{eqnarray}\label{Ricci}
\frac{\partial g}{\partial t}  = -2{\rm Ric} [g(t)],
\end{eqnarray}
with initial condition $g(0) = g_0$ where
$g(t) = c(t)  (\varphi _t)^* g_0$, and here $c(t)$ is a scaling parameter,
and $\varphi _t$ is a diffeomorphism.

Let $(N, g)$ be a nilpotent Lie group with left invariant metric $g$ and let $\mathfrak{n}$ be its conrredponded Lie
algebra with the induced inner product that we denote by $\langle.,.\rangle$. If not stated otherwise, we
assume that $N$ is 4-dimensional and that both metric $g$ and inner product $\langle.,.\rangle$ are of Lorentz
signature $(+,+,+,-)$.
%
\section{The Ricci soliton on $H_3 \times {\Bbb R}$ }
Algebra $\mathfrak{h}_3 \oplus {\Bbb R}$ is spanned by basis $\{x_1, x_2, x_3, x_4\}$ with nonzero commutator
\begin{align}
[x_1, x_2]=x_3
\end{align}
The algebra $\mathfrak{h}_3 \oplus {\Bbb R}$  2-step nilpotent with two dimensional center $Z(\mathfrak{h}_3 \oplus {\Bbb R} ) ={\cal L}(x_4, x_3)$
and one dimensional commutator subalgebra spanned by $x_3$.
Taking the group action on $H_3 \times {\Bbb R}$ for  coordinates $(x, y, z, w)$ is given by
\begin{eqnarray}
(x, y, z, w) \cdot (x', y', z', w')=(x+x', y+y'+xz', z+z', w+w')
\end{eqnarray}

The Lie algebra $\mathfrak{h}_3 \oplus {\Bbb R}$ of $H_3 \times {\Bbb R}$
has a basis with frame
\begin{eqnarray}\label{frame1}
{\bf X}_1  = \frac{\partial }{\partial x} ,\,\,\,\,\
{\bf X}_2  = \frac{\partial }{\partial y} +x\frac{\partial }{\partial z} ,\,\,\,\,\
{\bf X}_3  = \frac{\partial }{\partial z} ,\,\,\,\,\
{\bf X}_4  = \frac{\partial }{\partial w} ,
\end{eqnarray}
The coframe dual to the left invariant frame \eqref{frame1} is given by the basis of one-forms
\begin{eqnarray}\label{coframe1}
\omega_1  = dx ,\,\,\,\,\
\omega_2  = dy ,\,\,\,\,\
\omega_3  = dz-x dy ,\,\,\,\,\
\omega_4  = dw,
\end{eqnarray}
therefore 
\begin{eqnarray}\label{coframe102}
d\omega_1  = 0 ,\,\,\,\,\
d\omega_2  = 0 ,\,\,\,\,\
d\omega_3  = -\omega_1\wedge \omega_2 ,\,\,\,\,\
d\omega_4  = 0. 
\end{eqnarray}
\begin{thm}(\cite{Bokan2015})
Each left invariant Lorentz metric on the group $H_3 \times {\Bbb R}$, up to an automorphism
of $H_3 \times {\Bbb R}$, is isometric to one of the following
\begin{align*}
g_{\mu} = &\ dx^2 - dy^2 +\mu(x dy  - dz)^2+dw^2 , \\
g_{\lambda}^{\pm} = &\ dx^2 +dy^2 \pm\lambda(x dy  - dz)^2\mp dw^2, \\
g_0^1 = &\ dx^2 +dy^2 -2x dy\ dw +2dz\ dw, \\
g_0^2 = &\ dx^2 -2x\ dy^2 +dw^2 +2dy\ dz, \\
g_0^3 = &\ dx^2 +2dy\ dw +(x dy -dz)^2,
\end{align*}
\end{thm}
where $\lambda, \mu>0$.
It is easy to check
\begin{align*}
g_{\mu} = &\  \omega_1\otimes \omega_1 - \omega_2\otimes \omega_2 + \mu\ \omega_3\otimes \omega_3 +  \omega_4\otimes \omega_4 , \\
g_{\lambda}^{\pm} = &\  \omega_1\otimes \omega_1 + \omega_2\otimes \omega_2 \pm \lambda\ \omega_3\otimes \omega_3 \mp  \omega_4\otimes \omega_4 , \\
g_0^1  = &\  \omega_1\otimes \omega_1 +  \ \omega_2\otimes \omega_2 + 2 \omega_3\otimes \omega_4, \\
g_0^2 = &\  \omega_1\otimes \omega_1 + 2 \omega_2\otimes \omega_3 +  \omega_4\otimes \omega_4, \\
g_0^3  = &\  \omega_1\otimes \omega_1 + 2 \omega_2\otimes \omega_4 +  \omega_3\otimes \omega_3.
\end{align*}
We can write  $g_{\mu}$ and $g_{\lambda}^{\pm}$ metrics in general form 
\begin{eqnarray}\label{general1}
g= \omega_1\otimes \omega_1 + a_1\ \omega_2\otimes \omega_2 + a_2\ \omega_3\otimes \omega_3 + a_3\ \omega_4\otimes \omega_4,
\end{eqnarray}
where substituting $a_1=-a_3=-1,\; a_2=\mu$ in the \eqref{general1} we find $g_{\mu}$ and  substituting
 $a_1=1,\; a_2=\pm\lambda, \;a_3=\mp1$ in \eqref{general1}, leads to  $g_{\lambda}^{\pm}$.  
\begin{thm}
The left-invariant Lorentzian metrics  $g_{\mu}$ and $g_{\lambda}^{\pm}$ with general form \eqref{general1} satisfy a Ricci soliton equation
\begin{eqnarray}
2{\rm Ric}[g] + {\cal L}_{\bf X} g -\frac{3a_2}{a_1}g = 0,
\end{eqnarray}
where the vector field ${\bf X}$ is defined by
\begin{eqnarray}\nonumber
{\bf X} &=&\left( -\frac{a_2}{a_1}x+C_2y+C_3\right)\;{\bf X}_1-\left( \frac{1}{a_1}(C_2x+a_2y)+C_4\right)\;{\bf X}_2+\Big(\frac{C_2}{2a_1}(x^2 -a_1y^2)+C_3y+C_4x\\\label{vectfield1}
&&-\frac{a_2}{a_1}(xy-2z)+C_5\Big)\;{\bf X}_3+\left( \frac{a_2}{a_3}y-\frac{3a_2}{2a_1}w+C_1\right)\;{\bf X}_4.
\end{eqnarray}
Here for $g_{\mu}$ we have $a_1=-a_3=-1,\; a_2=\mu$ and  for $g_{\lambda}^{\pm}$
we have $a_1=1,\; a_2=\pm\lambda, \;a_3=\mp1$ and $C_i$s for $i=1, \ldots, 5$ are constants.
Therefore the left-invariant
Lorentzian metric  $g_{\lambda}^{+}$  is a shrinking Lorentz Ricci soliton while  $g_{\lambda}^{+}$  and  $g_{\mu}$ are expanding.
\end{thm}
{\noindent{\bf Proof:}}
For \eqref{general1}, there is a unique solution $g(t)$  to the Ricci flow
satisfied in \eqref{Ricci}, with the following form
\begin{eqnarray}\label{uniquesolution}
g(t)= f_1(t)\ \omega_1\otimes \omega_1 + f_2(t)\ \omega_2\otimes \omega_2 + f_3(t)\ \omega_3\otimes \omega_3 + f_4(t)\ \omega_4\otimes \omega_4,
\end{eqnarray}
satisfying the initial condition $g(0) = g_0$. The matrix of connection one-forms of a metric $g$ as in \eqref{general1} is
$$
\omega=(\omega^i_j)=\frac{1}{2}
\left(
  \begin{array}{cccc}
    0 & a_2\;\omega^3 & c\;\omega^2 & 0 \\[3mm]
    -\dfrac{a_2}{a_1}\;\omega^3 & 0 &  -\dfrac{a_2}{a_1}\;\omega^1 & 0 \\[3mm]
    -\omega^2 & \omega^1 & 0 & 0 \\[3mm]
   0 & 0 & 0 & 0 \\
  \end{array}
\right)
$$
and the matrix of curvature two-forms is given by
$$
\Omega=(\Omega^i_j)=\frac{1}{4}
\left(
  \begin{array}{cccc}
    0 & -3a_2\;\omega^1 \wedge \omega^2 &  \dfrac{a_2^2}{a_1}\;\omega^1 \wedge \omega^3 & 0 \\[3mm]
    \dfrac{3a_2}{a_1}\;\omega^1 \wedge \omega^2 & 0 &  \dfrac{a_2^2}{a_1}\;\omega^2 \wedge \omega^3 & 0 \\[3mm]
    -\dfrac{a_2}{a_1}\;\omega^1 \wedge \omega^3  & a_2\;\omega^2 \wedge \omega^3 & 0 & 0 \\[3mm]
   0 & 0 & 0 & 0 \\
  \end{array}
\right)
$$
therefore the Ricci curvature is
\begin{eqnarray}\label{Ricci1}
{\rm Ric}\;[g]= -\dfrac{a_2}{2a_1}\ \omega_1\otimes \omega_1  -\dfrac{a_2}{2}\ \omega_2\otimes \omega_2 + \dfrac{a_2^2}{2a_1}\ \omega_3\otimes \omega_3,
\end{eqnarray}
The scalar curvature is
 \begin{eqnarray}\label{scalarcurvature1}
{\rm S}\;[g]= -\dfrac{a_2}{2a_1},
\end{eqnarray}
and the corresponding Ricci operator is
 \begin{eqnarray}\label{Riccioperator1}
{\rm Rc}\;[g]= {\rm diag}\left(-\dfrac{3a_2}{8a_1}, -\dfrac{3a_2}{8}, \dfrac{5a_2^2}{8a_1},  \dfrac{a_2a_3}{8a_1}\right).
\end{eqnarray}

Now assume ${\bf X}$ be an arbitrary vector field on $H_3 \times {\Bbb R}$ by
\begin{eqnarray}\label{vectfield2}
{\bf X} = \sum_{i=1}^4 P^i(x, y, z, w)\;{\bf X}_i
\end{eqnarray}
 where the $P^i$s are smooth functions on $H_3 \times {\Bbb R}$. The component
functions of the Lie derivative ${\cal L}_{\bf X} g$ are recorded in the following symmetric matrix
$$
({\cal L}_{\bf X} g)_{ij}=
\left(
\begin{smallmatrix}
2P^1_x & P^1_y+(1+x^2)P^2_x-xP^3_x & P^3_x+P^1_z-xP^2_x & P^1_w+P^4_x \\[3mm]
   &~ 2xP^1-2xP^3_y+2(1+x^2)P^2_y ~ &~  P^3_y+(1+x^2)P^2_z-x(P^3_z+P^2_y)-P^1 ~&~ P^4_y-xP^3_w+(1+x^2)P^2_w \\[3mm]
    && 2P^3_z-2xP^2_z & P^4_z+P^3_w-xP^2_w \\[3mm] 
   & & & 2P^4_w
  \end{smallmatrix}
  \right)
$$ 
The corresponding Ricci soliton equation 
$2{\rm Ric}[g] + {\cal L}_{\bf X} g +\alpha g = 0$, concludes the following system of partial differential equations
\begin{eqnarray}\label{eq:determinsys}
\left\{ \begin{array}{lcl}
-\frac{a_2}{a_1}+({\cal L}_{\bf X} g)_{11}+\alpha=0,\\[2mm]
\frac{a_2^2}{a_1}x^2-a_2+({\cal L}_{\bf X} g)_{22}+a_1\alpha + a_2\alpha x^2=0,\\[2mm]
\frac{a_2^2}{a_1}+({\cal L}_{\bf X} g)_{33}+ a_2\alpha=0,\\[2mm]
0+({\cal L}_{\bf X} g)_{44}+a_3\alpha=0,\\[2mm]
-\frac{a_2^2}{a_1}x+({\cal L}_{\bf X} g)_{23}-a_2\alpha x=0,\\[2mm]
({\cal L}_{\bf X} g)_{ij}=0, \;\;\;\ other~ cases
\end{array}\right.
\end{eqnarray}
We can obtain the coefficient functions $P^i$s of the vector field \eqref{vectfield1}, by solving the above system of PDE's
and additionally it concludes $\alpha=-\dfrac{3a_2}{a_1}$. 
Putting $a_1=1,\; a_2=\lambda$ in metric  $g_{\lambda}^{+}$  we have $\alpha=-3\lambda$ and and then it is a shrinking Lorentz Ricci soliton while  in $g_{\lambda}^{-}$  and  $g_{\mu}$
substituting $a_1=1,\; a_2=-\lambda$  and $a_1=-1,\; a_2=\mu$ respectively, we have $\alpha=3\lambda$ and $\alpha=3\mu$ respectively that shows they
are expanding. 
%
\begin{thm}
The left-invariant Lorentzian metrics 
\begin{eqnarray}\label{eqg01}
g_0^1  =   \omega_1\otimes \omega_1 +  \ \omega_2\otimes \omega_2 + 2 \omega_3\otimes \omega_4,
\end{eqnarray}
 satisfies a Ricci soliton equation
\begin{eqnarray}\label{ricci1}
2{\rm Ric}[g_0^1] + {\cal L}_{\bf X} g_0^1 +\alpha g_0^1 = 0,
\end{eqnarray}
where the vector field ${\bf X}$ is defined by
\begin{eqnarray}\nonumber
{\bf X} &=&(C_1y-\frac{\alpha}{2} x+C_5\cos w+C_6\sin w +C_2)\;{\bf X}_1-(C_1x+\frac{\alpha}{2} y-C_5\sin w+C_6\cos w +C_3)\;{\bf X}_2\\\label{vectfield2}
&&+\left(\frac{C_1}{2}(x^2+y^2)+(C_3+\frac{\alpha y}{2})x +C_2y-\alpha z-\frac{w}{2}+C_4\right)\;{\bf X}_3+C_7\;{\bf X}_4,
\end{eqnarray}
where $C_i$s for $i=1, \ldots, 7$ are arbitrary constants and $\alpha\in {\Bbb R}$.
\end{thm}
{\noindent{\bf Proof:}}
 The matrix of connection one-forms of the metric $g_0^1$  is
$$
\omega=(\omega^i_j)=\frac{1}{2}
\left(
  \begin{array}{cccc}
    0 & \omega^4 & 0& \omega^2  \\[3mm]
    -\omega^4 & 0 &  0 & -\omega^1 \\[3mm]
    -\omega^2 & \omega^1 & 0 & 0 \\[3mm]
   0 & 0 & 0 & 0 \\
  \end{array}
\right)
$$
and the matrix of curvature two-forms is represented by
$$
\Omega=(\Omega^i_j)=\frac{1}{4}
\left(
  \begin{array}{cccc}
    0 & 0 &  0 & \omega^1 \wedge \omega^4 \\[3mm]
    0 & 0 & 0 & \omega^2 \wedge \omega^4\\[3mm]
    -\omega^1 \wedge \omega^4  & -\omega^2 \wedge \omega^4 & 0 & 0 \\[3mm]
   0 & 0 & 0 & 0 \\
  \end{array}
\right)
$$
Thus the Ricci curvature is
\begin{eqnarray}\label{Ricci1}
{\rm Ric}\;[g_0^1]= \dfrac{1}{2}\ \omega_4\otimes \omega_4,
\end{eqnarray}
and the corresponding scalar curvature is equal to zero. 
Therefore the Ricci operator is
 \begin{eqnarray}\label{Riccioperator1}
{\rm Rc}\;[g_0^1]= {\rm diag}\left(0, 0,  0, \dfrac{1}{2}\right).
\end{eqnarray}
Let ${\bf X}$ be an arbitrary vector field 
$\displaystyle{{\bf X} = \sum_{i=1}^4 P^i(x, y, z, w)\;{\bf X}_i}$
  on $H_3 \times {\Bbb R}$. 
The Lie derivative ${\cal L}_{\bf X} g_0^1$ as a symmetric matrix  is given by
$$
({\cal L}_{\bf X} g_0^1)_{ij}=
\left(
\begin{smallmatrix}
2P^1_x & P^1_y+P^2_x-xP^4_x&P^1_z+P^4_x & P^1_w+P^2+P^3_x \\[3mm]
   &~ 2P^2_y-2xP^4_y ~ &~  P^2_z+P^4_y-xP^4_z ~&~ -P^1+P^2_w+P^3_y-xP^4_w \\[3mm]
    && 2P^4_z & P^3_z+P^4_w \\[3mm]
   & & & 2P^3_w
  \end{smallmatrix}
  \right)
$$
Puttig in the Ricci soliton equation
\eqref{ricci1}, leads to
\begin{eqnarray}\label{eq:determinsys1}
\left\{ \begin{array}{lcl}
({\cal L}_{\bf X} g)_{11}+\alpha=0,\\[2mm]
({\cal L}_{\bf X} g)_{22}+\alpha =0,\\[2mm]
1+({\cal L}_{\bf X} g)_{44}=0,\\[2mm]
({\cal L}_{\bf X} g)_{34}+\alpha =0,\\[2mm]
({\cal L}_{\bf X} g)_{ij}=0, \;\;\;\ other~ cases
\end{array}\right.
\end{eqnarray}
Solving the  system \eqref{eq:determinsys1}, the coefficient functions $P^i$s of the vector field \eqref{vectfield2} are obtained.
%
\begin{thm}
The left-invariant Lorentzian metrics 
\begin{eqnarray}\label{eqg02}
g_0^2 = \omega_1\otimes \omega_1 + 2 \omega_2\otimes \omega_3 +  \omega_4\otimes \omega_4,
\end{eqnarray}
 satisfies a Ricci soliton equation
\begin{eqnarray}\label{riccig02}
2{\rm Ric}[g_0^2 ] + {\cal L}_{\bf X} g_0^2  +\alpha g_0^2  = 0,
\end{eqnarray}
where the vector field ${\bf X}$ is defined by
\begin{eqnarray}\nonumber
{\bf X} &=&\Big(-\frac{\alpha}{2}x-C_1 yw- C_2w+C_5z+\frac{C_5}{3}y^3+\frac{C_6}{2}y^2+C_7y+C_8\Big)\;{\bf X}_1+\Big( C_1w-C_5\big(\frac{y^2}{2}+x\big)\\\nonumber
&&-\frac{\alpha+2C_6}{4}y+C_9\Big)\;{\bf X}_2
+\Big(\big(\frac{C_6}{2}-\frac{3\alpha}{4}\big)z+\left(\frac{C_5}{12}y^2+\frac{C_6}{6}y+\frac{C_1}{2}w-\frac{C_5}{2}x+\frac{C_7}{2}\right)y^2+\big(\frac{\alpha}{4}-\frac{C_6}{2}\big)xy\\\nonumber
&&+C_8y+C_5yz-C_2yw-(C_7+C_9)x-C_3w+C_{10}\Big)\;{\bf X}_3+\Big(-\frac{\alpha}{2}w+C_1\big(xy-z+\frac{y^3}{6}\big)\\\label{vectfield3}
&&+C_2\big(x+\frac{y^2}{2}\big)+C_3y+C_4\Big)\;{\bf X}_4
\end{eqnarray}
Here  $C_i$s for $i=1, \ldots, 10$ and $\alpha$ are arbitrary constants.
\end{thm}
{\noindent{\bf Proof:}}
 Scince the matrix of connection one-forms of the metric $g_0^2$  is
$$
\omega=(\omega^i_j)=\frac{1}{2}
\left(
  \begin{array}{cccc}
    0 & \omega^2 & 0& 0  \\[3mm]
    0 & 0 &  0 & 0 \\[3mm]
    -\omega^2 & 0 & 0 & 0 \\[3mm]
   0 & 0 & 0 & 0 \\
  \end{array}
\right)
$$
then the matrix of curvature two-forms equals to zero matrix and  the Ricci curvature is
\begin{eqnarray}\label{Ricci2}
{\rm Ric}\;[g_0^2]= 0. 
\end{eqnarray}
Suppose that ${\bf X}$ be an arbitrary vector field 
$\displaystyle{{\bf X} = \sum_{i=1}^4 P^i(x, y, z, w)\;{\bf X}_i}$
  on $H_3 \times {\Bbb R}$. 
The Lie derivative ${\cal L}_{\bf X} g_0^2$ as a symmetric matrix  is given by
$$
({\cal L}_{\bf X} g_0^2)_{ij}=
\left(
\begin{smallmatrix}
2P^1_x & P^2+P^1_y-xP^2_x+P^3_x&P^1_z+P^2_x & P^1_w+P^4_x \\[3mm]
   & -2P^1+x2P^2_y+2P^3_y ~ &~  P^2_y-xP^2_z +P^3_z~&~ -xP^2_w+P^3_w+P^4_y \\[3mm]
    && 2P^2_z & P^2_w+P^4_z \\[3mm]
   & & & 2P^4_w
  \end{smallmatrix}
  \right)
$$
Puttig in the Ricci soliton equation
\eqref{riccig02}, leads to
\begin{eqnarray}\label{eq:determinsys3}
\left\{ \begin{array}{lcl}
({\cal L}_{\bf X} g)_{11}+\alpha=0,\\[2mm]
({\cal L}_{\bf X} g)_{44}+\alpha =0,\\[2mm]
({\cal L}_{\bf X} g)_{23}+\alpha =0,\\[2mm]
({\cal L}_{\bf X} g)_{ij}=0, \;\;\;\ other~ cases
\end{array}\right.
\end{eqnarray}
Solving  \eqref{eq:determinsys3}, one finds immediately the coefficient functions $P^i$s of the vector field \eqref{vectfield3}.
%
\begin{thm}
The left-invariant Lorentzian metrics 
\begin{eqnarray}\label{eqg03}
g_0^3  =  \omega_1\otimes \omega_1 + 2 \omega_2\otimes \omega_4 +  \omega_3\otimes \omega_3.
\end{eqnarray}
 satisfies a Ricci soliton equation
\begin{eqnarray}
2{\rm Ric}[g_0^3] + {\cal L}_{\bf X} g_0^3 +\alpha g_0^3 = 0,
\end{eqnarray}
where the vector field ${\bf X}$ is defined by
\begin{eqnarray}\nonumber
{\bf X} &=&\left( -\frac{\alpha}{2}x+\frac{C_1}{2}y^2+C_2y+C_3\right)\;{\bf X}_1+C_4\;{\bf X}_2+\Big(-\frac{\alpha}{2}z+\frac{C_1}{6}(y^3-6y)+\frac{C_2}{2}y^2+C_3y-C_4x\\\label{vectfield4}
&&+C_5\Big)\;{\bf X}_3+\left( \frac{y}{2}-\alpha w+\frac{C_1}{2}(2z-xy)-C_2x+C_6\right)\;{\bf X}_4.
\end{eqnarray}
$C_i$s for $i=1, \ldots, 6$  are arbitrary constants and $\alpha\in {\Bbb R}$.
\end{thm}
{\noindent{\bf Proof:}}
 The matrix of connection one-forms of the metric $g_0^3$  is
$$
\omega=(\omega^i_j)=\frac{1}{2}
\left(
  \begin{array}{cccc}
    0 & \omega^3 & \omega^2 & 0 \\[3mm]
   0 & 0 &  0 & 0 \\[3mm]
    -\omega^2 & \omega^1 & 0 & 0 \\[3mm] 
   -\omega^3 & 0 & -\omega^1 & 0 \\
  \end{array}
\right)
$$
and the matrix of curvature two-forms is represented by
$$
\Omega=(\Omega^i_j)=\frac{1}{4}
\left(
  \begin{array}{cccc}
    0 & -3\;\omega^1 \wedge \omega^2  &  0 & 0\\[3mm]
    0 & 0  & 0 & 0\\[3mm]
    0 & -\omega^2 \wedge \omega^3  & 0 & 0\\[3mm]
   3\;\omega^1 \wedge \omega^2 & 0 & \omega^2 \wedge \omega^3 & 0 \\
  \end{array}
\right)
$$
Thus the Ricci curvature tensor is
\begin{eqnarray}\label{Ricci1}
{\rm Ric}\;[g_0^1]= -\dfrac{1}{2}\ \omega_2\otimes \omega_2
\end{eqnarray}
and the corresponding scalar curvature  equals to zero. 
Therefore the Ricci operator is
 \begin{eqnarray}\label{Riccioperator1}
{\rm Rc}\;[g_0^1]= {\rm diag}\left(0, -\frac{1}{2},  0, 0\right).
\end{eqnarray}
Let ${\bf X}$ be an arbitrary vector field 
$\displaystyle{{\bf X} = \sum_{i=1}^4 P^i(x, y, z, w)\;{\bf X}_i}$
  on $H_3 \times {\Bbb R}$. 
The Lie derivative ${\cal L}_{\bf X} g_0^1$ as a symmetric matrix  is given by
$$
({\cal L}_{\bf X} g_0^3)_{ij}=
\left(
\begin{smallmatrix} 
2P^1_x~ & ~P^1_y-P^2_x-xP^3_x+P^4_x~ &~ P^1_z+P^2+P^3_x ~&~ P^1_w+P^2_x \\[3mm] 
  ~ &~ 2x(P^1-P^3_y)+2P^4_y  ~& ~ -P^1+P^3_y-xP^3_z+P^4_z ~&~ P^2_y-xP^3_w+P^4_w \\[3mm]
    &~ & 2P^3_z ~&~ P^2_z+P^3_w \\[3mm] 
   & ~& ~&~ 2P^2_w
  \end{smallmatrix} 
  \right)
$$
Puttig in the Ricci soliton equation
\eqref{ricci1}, leads to
\begin{eqnarray}\label{eq:determinsys2}
\left\{ \begin{array}{lcl}
({\cal L}_{\bf X} g)_{11}+\alpha=0,\\[2mm]
-1+({\cal L}_{\bf X} g)_{22}=0,\\[2mm]
({\cal L}_{\bf X} g)_{33}+\alpha=0,\\[2mm]
({\cal L}_{\bf X} g)_{24}+\alpha =0,\\[2mm]
({\cal L}_{\bf X} g)_{ij}=0, \;\;\;\ other~ cases
\end{array}\right.
\end{eqnarray}
The coefficient functions $P^i$s of the vector field \eqref{vectfield4}  are readily seen to be a solution to the system \eqref{eq:determinsys2}. 
\section{The Ricci soliton on $G_4$ }
The algebra ${\frak g}_4$  is spanned by basis $\{x_1, x_2, x_3, x_4\}$ with nonzero commutators
\begin{eqnarray}
[x_1, x_2] = x_3, \;\;\;\;\;\; [x_1, x_3] = x_4.
\end{eqnarray}
If ${\bf X}_1, {\bf X}_2, {\bf X}_3, {\bf X}_4$ are left invariant vector fields on $G_4$ defined by $x_1, x_2, x_3, x_4\in{\frak g}_4$, for
global coordinates $(x, y, z, w)$ on $G_4$ we have the relations
\begin{eqnarray}\label{frame2}
{\bf X}_1  = \frac{\partial }{\partial x} ,\,\,\,\,\
{\bf X}_2  = \frac{\partial }{\partial y} +x\frac{\partial }{\partial z}+\frac{x^2}{2} \frac{\partial }{\partial w} ,\,\,\,\,\
{\bf X}_3  = \frac{\partial }{\partial z}+x \frac{\partial }{\partial w} ,\,\,\,\,\
{\bf X}_4  = \frac{\partial }{\partial w} ,
\end{eqnarray}
The coframe dual to the left invariant frame \eqref{frame2} is obtaind by the basis of one-forms
\begin{eqnarray}\label{coframe12}
\omega_1  = dx ,\,\,\,\,\
\omega_2  = dy ,\,\,\,\,\
\omega_3  = dz-x\; dy ,\,\,\,\,\
\omega_4  = \frac{x^2}{2}\; dy-x\;dz+dw,
\end{eqnarray}
therefore 
\begin{eqnarray}\label{coframe12}
d\omega_1  =0 ,\,\,\,\,\
d\omega_2  = 0 ,\,\,\,\,\
d\omega_3  =-\omega_1\wedge \omega_2 ,\,\,\,\,\
d\omega_4  = -\omega_1\wedge \omega_3.
\end{eqnarray}
\begin{thm}(\cite{Bokan2015})
Each left invariant Lorentz metric on the group $G_4$, up to an automorphism
of $G_4$, is isometric to one of the following
\begin{align*}
g_{A}^\pm = &\pm dx^2 \mp dy^2 + a(xdy - dz)^2 - b(xdy - dz)(2dw - 2xdz + x^2dy) + \frac{c}{4}(2dw + x(xdy - 2dz))^2, \\
g_A = & dx^2 + dy^2 + a(xdy - dz)^2 - b(xdy - dz)(2dw - 2xdz + x^2dy)+\frac{c}{4}(2dw + x(xdy - 2dz))^2, \\
g_1^\lambda =& dx^2 + 2dwdy + xdy(xdy - 2dz) + \lambda(xdy - dz)^2, \\
g_2^\lambda =& 2dwdx + dy^2 + xdx(xdy - 2dz) + \lambda(xdy - dz)^2 , \\
g_3^\lambda =& dy^2 - 2dx(xdy - dz) + \frac{\lambda}{4}(2dw + x(xdy - 2dz))^2,\\
g_4^\lambda =& dx^2 - 2dy(xdy - dz) +  \frac{\lambda}{4}(2dw + x(xdy - 2dz))^2, 
\end{align*}
\end{thm}
where $\lambda>0$. With a simple calculation we can rewrite
\begin{align*}
g_{A}^\pm = &\ \pm \omega_1\otimes \omega_1 \mp \omega_2\otimes \omega_2 + a\ \omega_3\otimes \omega_3 +2b\ \omega_3\otimes \omega_4+ c\ \omega_4\otimes \omega_4 , \\
g_{A} = &\  \omega_1\otimes \omega_1 + \omega_2\otimes \omega_2 + a\ \omega_3\otimes \omega_3 +2b\ \omega_3\otimes \omega_4+ c\ \omega_4\otimes \omega_4 , \\
g_1^\lambda  = &\  \omega_1\otimes \omega_1 + 2 \ \omega_2\otimes \omega_4 + \lambda\ \omega_3\otimes \omega_3, \\
g_2^\lambda = &\  \omega_2\otimes \omega_2 + 2\ \omega_1\otimes \omega_4 + \lambda\ \omega_3\otimes \omega_3, \\
g_3^\lambda = &\  \omega_2\otimes \omega_2 + 2\ \omega_1\otimes \omega_3 + \lambda\ \omega_4\otimes \omega_4, \\
g_4^\lambda = &\  \omega_1\otimes \omega_1 + 2\ \omega_2\otimes \omega_3 + \lambda\ \omega_4\otimes \omega_4.
\end{align*}
%
\begin{thm}
The left-invariant Lorentzian metric $g_1^\lambda$ satisfies a Ricci soliton equation
\begin{eqnarray}
2{\rm Ric}[g_1^\lambda] + {\cal L}_{\bf X} g_1^\lambda +\alpha g_1^\lambda = 0,
\end{eqnarray}
where the vector field ${\bf X}$ is defined by
\begin{eqnarray}\label{vectfield22}
{\bf X} = -\frac{\alpha}{2}x\;{\bf X}_1+C_1\;{\bf X}_2-\Big(\frac{\alpha}{2}z+C_1x+C_2\Big)\;{\bf X}_3+\frac{1}{2}\left(\lambda y+\alpha (zx-2w)+C_1 x^2-2C_2 x+2C_3 \right)\;{\bf X}_4,
\end{eqnarray}
where $C_i$s for $i=1, \ldots, 3$   are arbitrary constants and $\alpha\in{\Bbb R}$.
\end{thm}
{\noindent{\bf Proof:}}
 The matrix of connection one-forms of a metric $g_1^\lambda$ is
$$
\omega=(\omega^i_j)=\frac{1}{2}
\left(
  \begin{array}{cccc}
    0 & (\lambda +1)\;\omega^3 & (\lambda +1)\;\omega^2 & 0 \\[3mm]
    0 & 0 &  0 & 0 \\[3mm]
    -\dfrac{1+\lambda}{\lambda}\;\omega^2 & \dfrac{\lambda-1}{\lambda}\;\omega^1 & 0 & 0 \\[4mm]
   -(1+\lambda)\;\omega^3  & 0 & (1-\lambda)\;\omega^1 & 0 \\
  \end{array}
\right)
$$
and the matrix of curvature two-forms is given by
$$
\Omega=(\Omega^i_j)=\frac{1}{4}
\left(
  \begin{array}{cccc}
    0 & -\dfrac{3\lambda^2+2\lambda-1}{\lambda}\;\omega^1 \wedge \omega^2 &  0 & 0 \\[5mm]
    0 & 0 &  0 & 0 \\[5mm]
    0  & -\dfrac{(1+\lambda)^2}{\lambda}\;\omega^2 \wedge \omega^3 & 0 & 0 \\[5mm]
   \dfrac{3\lambda^2+2\lambda-1}{\lambda}\;\omega^1 \wedge \omega^2 & 0 & (1+\lambda)^2\;\omega^2 \wedge \omega^3 & 0 \\
  \end{array}
\right)
$$
therefore the Ricci curvature is
\begin{eqnarray}\label{Ricci22}
{\rm Ric}\;[g_1^\lambda]= \dfrac{1-\lambda^2}{2\lambda}\ \omega_2\otimes \omega_2,
\end{eqnarray}
and the scalar curvature equals to zero. 
The corresponding Ricci operator is
 \begin{eqnarray}\label{Riccioperator1}
{\rm Rc}\;[g_1^\lambda]= {\rm diag}\left(0, \dfrac{1-\lambda^2}{2\lambda}, 0,0 \right).
\end{eqnarray}

Assume ${\bf X}$ be an arbitrary vector field on $G_4$ by
$\displaystyle {\bf X} = \sum_{i=1}^4 P^i(x, y, z, w)\;{\bf X}_i$
 where the $P^i$s are smooth functions on $G_4$. The component
functions of the Lie derivative ${\cal L}_{\bf X} g_1^\lambda$ are listed by
\begin{eqnarray*}
({\cal L}_{\bf X} g_1^\lambda)_{11}&=&2P^1_x,\\[2mm]
({\cal L}_{\bf X} g_1^\lambda)_{12}&=& P^1_y-\frac{3}{2}x^2P^2_x+x(P^2+P^3_x)(1-\lambda)+P^3+P^4_x, \\[2mm]
({\cal L}_{\bf X} g_1^\lambda)_{13}&=& P^1_z+\lambda P^2-2xP^2_x+\lambda P^3_x,\\ [2mm]
({\cal L}_{\bf X} g_1^\lambda)_{14}&=& P^1_w+P^2_x ,\\[2mm]
   ({\cal L}_{\bf X} g_1^\lambda)_{22}&=& 2(1+\lambda)xP^1 +x^2P^2_y-2\lambda x P^3_y+2P^4_y, \\[2mm]
  ({\cal L}_{\bf X} g_1^\lambda)_{23}&=& -(2+\lambda)P^1-2xP^2_y+\lambda P^3_y-\lambda xP^3_z-\frac{x^2}{2}P^2_z-P^3_z+P^4_z ,\\[2mm]
   ({\cal L}_{\bf X} g_1^\lambda)_{24}&=&   P^2_y-\frac{x^2}{2}P^2_w-(1+\lambda)xP^3_w+P^, \\[2mm]
    ({\cal L}_{\bf X} g_1^\lambda)_{33}&=& 2\lambda P^3_z-2xP_z^2, \\[2mm]
    ({\cal L}_{\bf X} g_1^\lambda)_{34}&=& P^2_z+\lambda P^3_w, \\[2mm]
   ({\cal L}_{\bf X} g_1^\lambda)_{44}&=& 2P^2_w.
\end{eqnarray*}
The corresponding Ricci soliton equation 
$2{\rm Ric}[g_1^\lambda] + {\cal L}_{\bf X} g_1^\lambda +\alpha g_1^\lambda = 0$, deduces the following system of PDE's
\begin{eqnarray}\label{eq:determinsys33}
\left\{ \begin{array}{lcl}
({\cal L}_{\bf X} g_1^\lambda)_{11}+\alpha=0,\\[2mm]
 \dfrac{1-\lambda^2}{\lambda}+({\cal L}_{\bf X} g_1^\lambda)_{22}=0,\\[2mm]
({\cal L}_{\bf X} g_1^\lambda)_{33}+ \lambda\alpha=0,\\[2mm]
({\cal L}_{\bf X} g)_{24}-a_2\alpha x=0,\\[2mm]
({\cal L}_{\bf X} g)_{ij}=0, \;\;\;\ other~ cases
\end{array}\right.
\end{eqnarray}
When we solve the system \eqref{eq:determinsys33},  the coefficient functions $P^i$s of the vector field \eqref{vectfield22} are found. 
%
\begin{thm}
The left-invariant Lorentzian metric $g_2^\lambda$ satisfies a Ricci soliton equation
\begin{eqnarray}
2{\rm Ric}[g_2^\lambda] + {\cal L}_{\bf X} g_2^\lambda +\alpha g_2^\lambda = 0,
\end{eqnarray}
for $\alpha=0$ where the vector field ${\bf X}$ is defined by
\begin{eqnarray}\nonumber
{\bf X} &=&- \left(\frac{C_1}{2}\big(x^2+\frac{2}{\lambda}\big)+C_2x+C_3\right)\;{\bf X}_2+\Big(\frac{C_1}{6}x^3+\frac{C_2}{2}x^2+C_3x+C_4\Big)\;{\bf X}_3+
\Big(-\frac{C_1}{24}x^4-\frac{C_2}{2}x^3\\\label{vectfield33}
&&-\frac{C_3}{2}x^2+(\frac{\lambda}{2}-C_4)x+C_2y+C_1z+C_5\Big)\;{\bf X}_4,
\end{eqnarray}
where $C_i$s for $i=1, \ldots, 5$   are arbitrary constants and then $g_2^\lambda$ is a shrinking. 
\end{thm}
{\noindent{\bf Proof:}}
 The matrix of connection one-forms of a metric $g_2^\lambda$ is  
$$
\omega=(\omega^i_j)=\frac{1}{2}
\left(
  \begin{array}{cccc}
    0 &0& 0 & 0 \\[3mm]
    -\lambda\ \omega^3  & 0 &   -\lambda\ \omega^1 & 0 \\[3mm]
    -\dfrac{1}{\lambda}\ \omega^1- \omega^2 & -\omega^2 & 0 & 0 \\[4mm] 
   0  & \lambda\ \omega^3 & \omega^1+\lambda \omega^2 & 0 \\
  \end{array}
\right)
$$
and the matrix of curvature two-forms is given by
$$
\Omega=(\Omega^i_j)=\frac{1}{4}
\left(
  \begin{array}{cccc}
    0 & 0 &  0 & 0 \\[5mm]
    3\lambda\;\omega^1 \wedge \omega^2 & 0 &  0 & 0 \\[5mm]
    -\lambda\;\omega^1 \wedge \omega^3  & 0 & 0 & 0 \\[5mm]
    0& -3\lambda\;\omega^1 \wedge \omega^2  & \lambda^2\;\omega^1 \wedge \omega^3 & 0 \\
  \end{array}
\right),
$$
therefore the Ricci curvature is 
\begin{eqnarray}\label{Ricci22} 
{\rm Ric}\;[g_2^\lambda]= -\dfrac{\lambda}{2}\ \omega_1\otimes \omega_1,
\end{eqnarray}
and the scalar curvature is equal to zero. 
The corresponding Ricci operator is
 \begin{eqnarray}\label{Riccioperator22} 
{\rm Rc}\;[g_2^\lambda]= {\rm diag}\left(-\dfrac{\lambda}{2}, 0, 0,0 \right). 
\end{eqnarray}

Now assume ${\bf X}$ be an arbitrary vector field on $G_4$ by
$\displaystyle {\bf X} = \sum_{i=1}^4 P^i(x, y, z, w)\;{\bf X}_i$
 where the $P^i$s are smooth functions on $G_4$. The component
functions of the Lie derivative symmetric matrix ${\cal L}_{\bf X} g_2^\lambda$ are listed by
\begin{eqnarray*}
({\cal L}_{\bf X} g_2^\lambda)_{11}&=&2P^4_x+2P^3,\\
 ({\cal L}_{\bf X} g_2^\lambda)_{12}&=&P^4_y-xP^3_y-\frac{x^2}{2}P^2_y-\lambda xP^3_x+x^2P^1_x-\lambda xP^2+P^2_x+2xP^1,  \\
({\cal L}_{\bf X} g_2^\lambda)_{13}&=&P^4_z-xP^3_z-\frac{x^2}{2}P^2_z+\lambda P^3_x-2xP^1_x+\lambda P^2-2P^1,\\
({\cal L}_{\bf X} g_2^\lambda)_{14}&=&P^4_w-xP^3_w-\frac{x^2}{2}P^2_w+P^1_x,\\
({\cal L}_{\bf X} g_2^\lambda)_{22}&=&-2\lambda xP^3_y+x^2P^1_y+2\lambda xP^1+2P^2_y,\\
({\cal L}_{\bf X} g_2^\lambda)_{23}&=&-\lambda xP^3_z+\lambda P^3_y-2xP^1_y-\lambda P^1+P^2_z,\\
({\cal L}_{\bf X} g_2^\lambda)_{24}&=&P^1_y+P^2_w-\lambda xP^3_w,\\
({\cal L}_{\bf X} g_2^\lambda)_{33}&=&-2xP^1_z+2\lambda P^3_z,\\
({\cal L}_{\bf X} g_2^\lambda)_{34}&=&P^1_z+\lambda P^3_w,\\
({\cal L}_{\bf X} g_2^\lambda)_{44}&=&2P^1_w,
 \end{eqnarray*} 
The corresponding Ricci soliton equation with $\alpha=0$,  is 
$2{\rm Ric}[g_2^\lambda] + {\cal L}_{\bf X} g_2^\lambda  = 0$, and it leads to the following system 
\begin{eqnarray}\label{eq:determinsys44}
\left\{ \begin{array}{lcl}
-\lambda+({\cal L}_{\bf X} g_2^\lambda)_{44}=0,\\[2mm] 
({\cal L}_{\bf X} g_2^\lambda)_{ij}=0. \;\;\;\ other~ cases
\end{array}\right.
\end{eqnarray}
Solving the above system \eqref{eq:determinsys44},  we can find the vector field \eqref{vectfield33}. 
\section{Acknowledgement}
I wish to express my  sincere gratitude to Dr. Thomas H. Wears from Longwood University  for 
some helpful discussion about the paper and sending her student thesis for me. 
I Also acknowledge  the funding support of  Babol Noshirvani University of Technology under
Grant No. BNUT/391024/1400.

%
\end{document}